# New identities with Stirling, hyperharmonic, and derangement numbers, Bernoulli and Euler polynomials, powers, and factorials

Khristo N. Boyadzhiev

*Department of Mathematics*
*Ohio Northern University*
*Ada, OH 45810, USA*

[k-boyadzhiev@onu.edu](k-boyadzhiev@onu.edu)

**Abstract**

A number of identities are proved by using Stirling transforms. These identities involve Stirling numbers of the first and second kinds, hyperharmonic and derangement numbers, Bernoulli and Euler numbers and polynomials, powers, power sums, and factorials.

**1. The Stirling transform**

Throughout $s(n,k)$ are the Stirling numbers of the first kind and $S(n,k)$ are the Stirling numbers of the second kind. These numbers are very popular in combinatorics, as demonstrated by Comtet in [11], Graham et al in [13], and numerous other publications. At the same time they have important applications in analysis due to their connection to Newton series, asymptotic expansions, and series summation (see [5, 6, 7, 10, 14, 16, 18]). It is beneficial to study their connection to other important numbers such as Bernoulli numbers, Euler numbers, derangement numbers, and harmonic numbers. One way to find such connections is by using the Stirling transform.

Let $\{a_n\}_{n=0}^{\infty}$ be a number sequence. Its Stirling transform is another sequence defined by

$$b_n = \sum_{k=0}^{n} S(n,k) a_k$$

The inversion formula is

$$a_n = \sum_{k=0}^{n} s(n,k) b_k .$$

The inversion is based on the identities

$$\sum_{k=0}^{n} S(n,k) s(k,j) = \begin{cases} 1 & n = j \\ 0 & n \neq j \end{cases} \quad \text{and} \quad \sum_{k=0}^{n} s(n,k) S(k,j) = \begin{cases} 1 & n = j \\ 0 & n \neq j \end{cases}.$$

(see [11], [13], [16].).

An interesting combinatorial interpretation of this transform was given by Bernstein and Sloane in [2]. Rahmani [17] considered a generalization based on the $r$-Stirling numbers. The recent book [4] contains a short table of Stirling transform formulas.

The Stirling transform is related to a couple of interesting series transformations. Namely, if

(1) $$f(t) = \sum_{n=0}^{\infty} \frac{a_n}{n!} t^n$$

is an exponential generating function we have the series transformation formulas (Comtet [11])

(2) $$f\left(\frac{\mu}{\lambda}(e^{\lambda t} - 1)\right) = \sum_{n=0}^{\infty} \frac{t^n}{n!} \left\{ \sum_{k=0}^{n} S(n,k) \lambda^{n-k} \mu^k a_k \right\}$$

(3) $$f\left(\frac{\mu}{\lambda} \log(1 + \lambda t)\right) = \sum_{n=0}^{\infty} \frac{t^n}{n!} \left\{ \sum_{k=0}^{n} s(n,k) \lambda^{n-k} \mu^k a_k \right\}$$

where $\lambda$ and $\mu$ are parameters.

Interesting identities in the form of Stirling transform were obtained by Merlini et al in [15] for the Cauchy numbers introduced in [11]. In this paper we present several Stirling transform identities involving Euler and Bernoulli polynomials and numbers, derangement numbers, hyperharmonic numbers, powers and factorials. In particular, we evaluate the expressions

$$\sum_{k=0}^{n} S(n,k)(-1)^k k! h_k^{(p)}, \quad \sum_{k=0}^{n} s(n,k) B_k(x), \quad \sum_{k=1}^{n} s(n,k) B_{k-1}(-1)^k, \quad \sum_{k=0}^{n} s(n,k) E_k(x)$$

where $h_k^{(p)}$ are the hyperharmonic numbers, $B_k(x)$ are Bernoulli's polynomials, $B_k$ are the Bernoulli numbers, and $E_n(x)$ are Euler's polynomials. The first sum for $p=1$ becomes

$$\sum_{k=0}^{n} S(n,k)(-1)^k k! H_k = (-1)^n n$$

where $H_n$ are the harmonic numbers. As corollaries we obtain a representation of $h_k^{(p)}$, $B_n(x)$, and $E_n(x)$ in terms of the Stirling numbers $S(n,k)$. In Proposition 7 we show a connection among Stirling numbers, Bernoulli numbers, and harmonic numbers

$$\sum_{k=1}^{n} s(n,k) B_{k-1} = (-1)^{n-1}(n-1)! H_n.$$

Further, we evaluate also the sums

$$\sum_{k=0}^{n} S(n,k)(-1)^k D_k, \quad \sum_{k=1}^{n} S(n,k)\frac{x^k}{k}, \quad \sum_{k=1}^{n} S(n,k)\frac{x^k}{k^2}, \quad \sum_{k=1}^{n} S(n,k)(k-1)! x^k$$

where $D_k$ are the derangement numbers. We also give a recurrence formula for the sums $\sum_{k=1}^{n} S(n,k) k^p$ for $p=1,2,...$. In Section 5 we prove the identity

$$\sum_{k=1}^{p+1} S(p+1,k)\frac{x^k}{k} = e^{-x} \sum_{n=0}^{\infty} \left(1^p + 2^p + ... + n^p\right)\frac{x^n}{n!}$$

connecting Stirling numbers and sums of powers. At the end we also compute the alternating binomial transform of the exponential polynomials $\varphi_k(x)$

$$\sum_{k=0}^{n}\binom{n}{k}(-1)^{k}\varphi_{k}(x)=1+x\sum_{j=0}^{n-1}(-1)^{j+1}\varphi_{j}(x).$$

## 2. Stirling transform of hyperharmonic numbers

Given the harmonic numbers $H_n = 1 + \frac{1}{2} + \ldots + \frac{1}{n}$, the hyperharmonic numbers $h_n^{(p)}$ are defined in the following way: Set $h_n^{(1)} = H_n$ and further,

$$h_n^{(2)} = H_1 + H_2 + \ldots + H_n, \quad h_n^{(3)} = h_1^{(2)} + h_2^{(2)} + \ldots + h_n^{(2)}, \ldots$$

etc. (see Conway and Guy [12]). For convenience we also set $h_0^{(p)} = 0$. The hyperharmonic numbers can be written in the compact form

$$h_n^{(p+1)} = \binom{n+p}{n}(H_{n+p} - H_p), \ (p = 0, 1, 2, \ldots).$$

**Theorem 1.** *For $p = 0, 1, \ldots$ we have*

(4) $$\sum_{k=0}^{n} S(n,k)(-1)^k k! h_k^{(p)} = (-1)^n n p^{n-1}$$

*In particular, when $p = 1$ we have $h_n^{(1)} = H_n$ and*

(5) $$\sum_{k=0}^{n} S(n,k)(-1)^k k! H_k = (-1)^n n.$$

*Proof.* The generating function for the hyperharmonic numbers is given by

$$\frac{-\ln(1-t)}{(1-t)^p} = \sum_{n=0}^{\infty} h_n^{(p)} t^n$$

for $|t| < 1$ (see equation (7.43) in [13]). We will write this in the form

(6) $$\frac{-\ln(1+t)}{(1+t)^p} = \sum_{n=0}^{\infty} (-1)^n n! h_n^{(p)} \frac{t^n}{n!}$$

to match the structure of (1). Now we apply property (2) with $\lambda = \mu = 1$ to the function $-\ln(1+t)/(1+t)^p$ in (6) replacing there $t$ by $e^t - 1$. The result is

$$\frac{-\ln(1+(e^t-1))}{(1+(e^t-1))^p} = -te^{-pt} = \sum_{m=0}^{\infty} \frac{(-1)^{m+1} p^m t^{m+1}}{m!}$$

$$= \sum_{n=1}^{\infty} (-1)^n n p^{n-1} \frac{t^n}{n!} = \sum_{n=0}^{\infty} \frac{t^n}{n!} \left\{ \sum_{k=0}^{n} S(n,k)(-1)^k k! h_k^{(p)} \right\}.$$

Comparing coefficients in the last two sums gives identity (4) of the theorem. Note that the first term in the last sum is zero, since $h_0^{(p)} = 0$. With this the theorem is proved.

The identity in equation (5) was obtained previously by Kellner in [23, Lemma 1.1].

**Corollary 2.** *By inversion of the Stirling transform we find from (4) and (5) also*

(7) $\quad \sum_{k=0}^{n} s(n,k)(-1)^k k p^{k-1} = (-1)^n n! h_n^{(p)}$, or

$$h_n^{(p)} = \frac{1}{n!} \sum_{k=0}^{n} s(n,k)(-1)^{n-k} k p^{k-1}$$

*for any integer $p \geq 0$. In particular, when $p = 1$,*

$$\sum_{k=0}^{n} s(n,k)(-1)^k k = (-1)^n n! H_n.$$

## 3. Stirling transform of Bernoulli and Euler polynomials

The Euler polynomials are defined by the generating function

(8) $\quad \dfrac{2e^{xt}}{e^t+1} = \sum_{n=0}^{\infty} E_n(x) \dfrac{t^n}{n!}$.

We shall see that $E_n(x)$ can be written in terms of $S(n,k)$.

**Theorem 3.** *For any* $n \geq 0$

$$\sum_{k=0}^{n} s(n,k) E_k(x) = n! \sum_{k=0}^{n} \binom{x}{k} \left(-\frac{1}{2}\right)^{n-k}$$

$$E_n(x) = \sum_{k=0}^{n} S(n,k) k! \left\{ \sum_{j=0}^{k} \binom{x}{j} \left(-\frac{1}{2}\right)^{k-j} \right\}.$$

The following well-known technical result will be needed.

**Lemma 4.** *Let* $g(t) = a_0 + a_1 t + \ldots$ *be a power series. Then*

$$\frac{g(t)}{1-\lambda t} = \sum_{n=0}^{\infty} t^n \left\{ \sum_{k=0}^{n} a_k \lambda^{n-k} \right\} \quad \text{and} \quad \frac{g(t)}{1+\lambda t} = \sum_{n=0}^{\infty} t^n \left\{ \sum_{k=0}^{n} a_k (-1)^{n-k} \lambda^{n-k} \right\}$$

*where* $\lambda$ *is a parameter.*

(For the proof we expand $(1-\lambda t)^{-1}$ as a geometric power series for $|\lambda t| < 1$ and multiply the two power series. Then we use analytic continuation to drop the restriction.)

*Proof of the theorem.* In the series transformation formula (3) we set $\lambda = \mu = -1$. The formula becomes

$$f(\log(1-t)) = \sum_{n=0}^{\infty} \frac{t^n}{n!} \left\{ (-1)^n \sum_{k=0}^{n} s(n,k) a_k \right\}$$

and we apply it to the generating function (8) of $E_n(x)$, replacing $t$ by $\ln(1-t)$ in $2e^{xt}/(e^t+1)$. This gives the representation

$$\frac{2(1-t)^x}{2-t} = \frac{1}{1-t/2}(1-t)^x = \frac{1}{1-t/2} \sum_{n=0}^{\infty} \binom{x}{n} (-1)^n t^n = \sum_{n=0}^{\infty} \frac{t^n}{n!} \left\{ (-1)^n \sum_{k=0}^{n} s(n,k) E_k(x) \right\}.$$

According to Lemma 4 with $\lambda = 1/2$ the first series becomes

$$\frac{1}{1-t/2} \sum_{n=0}^{\infty} \binom{x}{n} (-1)^n t^n = \sum_{n=0}^{\infty} t^n \left\{ \sum_{k=0}^{n} \binom{x}{k} \frac{(-1)^k}{2^{n-k}} \right\} = \sum_{n=0}^{\infty} \frac{t^n}{n!} \left\{ n! \sum_{k=0}^{n} \binom{x}{k} \frac{(-1)^k}{2^{n-k}} \right\}.$$

Comparing coefficients we obtain the first formula of the proposition. The second one follows by inversion of the Stirling transform. The proposition is proved.

The Euler numbers $E_n$ can be defined as $E_n = E_n\left(\frac{1}{2}\right)$. Therefore, we have the representation of these numbers in terms of Stirling numbers for every $n \geq 0$

$$E_n = \sum_{k=0}^{n} S(n,k) k! \left\{ \sum_{j=0}^{k} \binom{1/2}{j} \left(-\frac{1}{2}\right)^{k-j} \right\}.$$

In view of the identity

$$\binom{1/2}{j} = \binom{2j}{j} \frac{(-1)^{j+1}}{2^{2j}(2j-1)}$$

we can write the above representation in the form

$$(9) \qquad E_n = \sum_{k=0}^{n} S(n,k) k! (-1)^k \left\{ \sum_{j=0}^{k} \binom{2j}{j} \frac{1}{2^{k+j}(1-2j)} \right\}.$$

This formula connects together three important sets of numbers – the Euler numbers, the Stirling numbers, and the central binomial coefficients $\binom{2j}{j}$.

Recently Boutiche, Rarmani, and Srivastava obtained similar representations for the generalized Euler polynomials $E_n^{(\alpha)}(x)$ (see Theorem 3.1 in [3]). The representation in [3] is different and coincides with ours only when $\alpha = 1$ and $x = 0$. In another recent paper [19] the same authors proved representations for Frobenius-Euler polynomials of higher order.

A similar representation is true for the Bernoulli polynomials. The Bernoulli polynomials can be defined through their generating function

$$(10) \qquad \frac{te^{xt}}{e^t - 1} = \sum_{n=0}^{\infty} B_n(x) \frac{t^n}{n!}$$

and $B_n = B_n(0)$ are the Bernoulli numbers [18].

It is well known that the Bernoulli numbers can be expressed in terms of Stirling numbers of the second kind (see [9] for comments)

$$B_n = \sum_{k=0}^{n} S(n,k) \frac{(-1)^k k!}{k+1}$$

and using the representation

$$B_n(x) = \sum_{p=0}^{n} \binom{n}{p} B_p x^{n-p}$$

we can write after changing the order of summation

$$B_n(x) = \sum_{k=0}^{n} \frac{(-1)^k k!}{k+1} \left\{ \sum_{p=0}^{n} \binom{n}{p} S(p,k) x^{n-p} \right\}$$

that is, expressing the Bernoulli polynomials in terms of Stirling numbers. However, the polynomials $B_n(x)$ can be expressed in terms of $S(n,k)$ as a Stirling transform formula. The following result was obtained by Todorov [21], equation (49), in a somewhat different form. Later this was extended by Rahmani [17] (the formula there, however, has a missing factor). We shall give a new proof here for completeness.

**Theorem 5**. *For any* $n \geq 0$

$$\sum_{k=0}^{n} s(n,k) B_k(x) = n! \sum_{k=0}^{n} \binom{x}{k} \frac{(-1)^{n-k}}{n-k+1}$$

$$B_n(x) = \sum_{k=0}^{n} S(n,k) k! \left\{ \sum_{j=0}^{k} \binom{x}{j} \frac{(-1)^{k-j}}{k-j+1} \right\}.$$

When $x = 0$ we have

$$B_n = B_n(0) = \sum_{k=0}^{n} S(n,k) k! \frac{(-1)^k}{k+1}$$

which is the well-known representation mentioned above.

For the proof of this theorem we need a technical result. Integrating with respect to $\lambda$ the second equation in Lemma 4 and setting $\lambda = 1$ we find

$$\frac{\ln(1+t)}{t}\sum_{n=0}^{\infty} a_n t^n = \sum_{n=0}^{\infty} t^n \left\{ \sum_{k=0}^{n} \frac{a_k (-1)^{n-k}}{n-k+1} \right\}$$

for any power series $a_0 + a_1 t + a_2 t^2 + \dots$.

*Proof of the theorem.* We apply property (3) with $\lambda = \mu = 1$ to the generating function (10), replacing $t$ by $\ln(1+t)$. The result is, according to the above corollary

$$\frac{\ln(1+t)(1+t)^x}{t} = \frac{\ln(1+t)}{t} \sum_{n=0}^{\infty} \binom{x}{n} t^n = \sum_{n=0}^{\infty} t^n \left\{ \sum_{k=0}^{n} \binom{x}{k} \frac{(-1)^{n-k}}{n-k+1} \right\}$$

$$= \sum_{n=0}^{\infty} \frac{t^n}{n!} \left\{ n! \sum_{k=0}^{n} \binom{x}{k} \frac{(-1)^{n-k}}{n-k+1} \right\} = \sum_{n=0}^{\infty} \frac{t^n}{n!} \left\{ \sum_{k=0}^{n} s(n,k) B_k(x) \right\}.$$

Comparing coefficients we come to the first equation in the theorem. The second equation follows from this by inversion of the Stirling transform. The proof is complete.

The next result is even more interesting, because it connects together three sets of special numbers: namely, Stirling, Bernoulli, and harmonic numbers.

**Theorem 6**. *For every $n \geq 1$,*

$$\sum_{k=1}^{n} s(n,k) B_{k-1} = (-1)^{n-1} (n-1)! H_n$$

$$B_{n-1} = \sum_{k=1}^{n} S(n,k)(-1)^{k-1}(k-1)! H_k \quad \text{or} \quad B_n = \sum_{k=1}^{n+1} S(n+1,k)(-1)^{k-1}(k-1)! H_k$$

$$\sum_{k=1}^{n} s(n,k) B_{k-1} (-1)^k = \frac{(-1)^n n!}{n^2}$$

$$B_{n-1} = (-1)^n \sum_{k=1}^{n} S(n,k) \frac{k!}{k^2}(-1)^k \quad \text{or} \quad B_n = (-1)^{n+1} \sum_{k=1}^{n+1} S(n+1,k) \frac{k!}{k^2}(-1)^k.$$

*Proof.* We set $x = 0$ in (10) and integrate to obtain the representation

$$\text{Li}_2(1-e^{-t}) = \sum_{m=0}^{\infty} \frac{B_m t^{m+1}}{m!(m+1)} = \sum_{n=1}^{\infty} B_{n-1} \frac{t^n}{n!}$$

where $\text{Li}_2(x) = \sum_{n=1}^{\infty} x^n / n^2$ is the dilogarithm. We now apply property (3) to the function $f(t) = \text{Li}_2(1-e^{-t})$ taking $\lambda = \mu = -1$. That is, we replace $t$ by $\ln(1-t)$ to get

$$\text{Li}_2\left(\frac{-t}{1-t}\right) = \sum_{n=1}^{\infty} \frac{t^n}{n!} \left\{ (-1)^n \sum_{k=1}^{n} s(n,k) B_{k-1} \right\}.$$

At the same time it easy to see that

$$\text{Li}_2\left(\frac{-t}{1-t}\right) = -\sum_{n=1}^{\infty} H_n \frac{t^n}{n} = -\sum_{n=1}^{\infty} (n-1)! H_n \frac{t^n}{n!}.$$

(A hint for the proof: Differentiating this equation brings to the generating function of the harmonic numbers $\frac{-\ln(1-t)}{1-t} = \sum_{n=1}^{\infty} H_n t^n$.) The first equation in the proposition follows from here by comparing coefficients. The second equation is just the inversion of the first one.

For the third equation we take $\lambda = 1$, $\mu = -1$ in (3) and replace $t$ by $-\ln(1+t)$ in $\text{Li}_2(1-e^{-t})$ to get

$$\text{Li}_2(-t) = \sum_{n=1}^{\infty} \frac{(-1)^n n! \, t^n}{n^2 \, n!} = \sum_{n=1}^{\infty} \frac{t^n}{n!} \left\{ \sum_{k=1}^{n} s(n,k) B_{k-1}(-1)^k \right\}$$

and the third equation follows by comparing coefficients. The last equation is again by inversion.

The theorem is proved.

(See also [23, Proposition 2.2] for similar results.)

## 4. Stirling transforms of powers

In this section we provide formulas for evaluating the numbers

$$M(n, p) = \sum_{k=0}^{n} S(n,k) k^p$$

for $p = 0, 1, \ldots$. We give a recurrence relation for these numbers from which they can be computed in a consecutive manner. The evaluation is given in terms of the Bell numbers

$$b_n = \sum_{k=0}^{n} S(n,k).$$

In what follows we shall need use the exponential polynomials [6], [7].

$$\varphi_n(x) = \sum_{k=0}^{n} S(n,k) x^k.$$

and some of their properties. In particular, with the notation $D = \dfrac{d}{dx}$ we have

(11) $\quad \varphi_p(x) = e^{-x}(xD)^p e^x = e^{-x} \sum_{n=0}^{\infty} \dfrac{n^p x^n}{n!}$

(12) $\quad \varphi_{p+1}(x) = xD\varphi_p(x) + x\varphi_p(x).$

It is clear that $b_n = \varphi_n(1)$ and

$$M(n, p) = (xD)^p \varphi_n(x) \big|_{x=1}.$$

**Theorem 7.** *For any $n \geq 0$, $p \geq 0$*

(13) $\quad M(n, p+1) = M(n+1, p) - \sum_{j=0}^{p} \binom{p}{j} M(n, j).$

*We have $M(n, 0) = b_n$. Using the recurrence (13) we compute*

$$M(n, 1) = \sum_{k=0}^{n} S(n,k) k = b_{n+1} - b_n$$

$$M(n, 2) = \sum_{k=0}^{n} S(n,k) k^2 = b_{n+2} - 2b_{n+1}$$

$$M(n, 3) = \sum_{k=0}^{n} S(n,k) k^3 = b_{n+3} - 3b_{n+2} + b_n$$

$$M(n,4) = \sum_{k=0}^{n} S(n,k)k^4 = b_{n+4} - 4b_{n+3} + 4b_{n+1} + b_n$$

$$M(n,5) = \sum_{k=0}^{n} S(n,k)k^5 = b_{n+5} - 5b_{n+4} + 10b_{n+2} + 5b_{n+1} - 2b_n$$

*etc.*

For the proof we need the following lemma which helps to evaluate $(xD)^p \varphi_n$ by recurrence.

**Lemma 8**. *For every $n \geq 0$ and $p \geq 0$*

(14) $\quad (xD)^{p+1} \varphi_n(x) = (xD)^p \varphi_{n+1}(x) - x \sum_{j=0}^{p} \binom{p}{j} (xD)^j \varphi_n(x).$

*Proof.* The simple proof follows from property (12) and the Leibnitz formula for the differential operator $xD$. We have

$$(xD)^{p+1} \varphi_n(x) = (xD)^p \varphi_{n+1}(x) - (xD)^p (x \varphi_n(x))$$

and

$$(xD)^p (x \varphi_n(x)) = \sum_{j=0}^{p} \binom{p}{j} \{(xD)^{p-j} x\} \{(xD)^j \varphi_n(x)\} = x \sum_{j=0}^{p} \binom{p}{j} \{(xD)^j \varphi_n(x)\}$$

since $(xD)^j x = x$ for every nonnegative integer $j$. The proof is finished.

Setting $x = 1$ in the lemma proves Theorem 7.

Applying the lemma for $p = 0, 1, ...$ and using the property $(xD) \varphi_n(x) = \varphi_{n+1}(x) - x \varphi_n(x)$ we find consecutively

(15) $\quad \sum_{k=0}^{n} S(n,k) k x^k = (xD) \varphi_n(x) = \varphi_{n+1}(x) - x \varphi_n(x)$

$$\sum_{k=0}^{n} S(n,k) k^2 x^k = (xD)^2 \varphi_n(x) = \varphi_{n+2}(x) - 2\varphi_{n+1}(x) + (x^2 - x) \varphi_n(x)$$

*etc.*

## 5. Transforms of reciprocals and connection to power sums

**Proposition 9**. *For every $p \geq 0$*

(16) $\quad \sum_{k=1}^{p+1} S(p+1,k) \frac{x^k}{k} = e^{-x} \sum_{n=0}^{\infty} \left(1^p + 2^p + ... + n^p\right) \frac{x^n}{n!}.$

*Proof.* We shall see that both sides in this equation have the same derivative. Then since both sides are zeros for $x=0$, they are equal. Differentiating both sides we find

$$\sum_{k=1}^{p+1} S(p+1,k)x^{k-1} = -e^{-x}\sum_{n=0}^{\infty}\left(1^p + 2^p + \ldots + n^p\right)\frac{x^n}{n!}$$

$$+e^{-x}\sum_{n=1}^{\infty}\left(1^p + 2^p + \ldots + n^p\right)\frac{x^{n-1}}{(n-1)!} = e^{-x}\sum_{m=0}^{\infty}(m+1)^p\frac{x^m}{m!}.$$

(after setting $m = n-1$ in the second sun). Further, we multiply both sides by $x$ to get

$$\sum_{k=1}^{p+1} S(p+1,k)x^k = e^{-x}\sum_{m=0}^{\infty}(m+1)^p\frac{x^{m+1}}{m!} = e^{-x}\sum_{m=0}^{\infty}(m+1)^{p+1}\frac{x^{m+1}}{(m+1)!}$$

$$= e^{-x}\sum_{m=0}^{\infty}(m+1)^{p+1}\frac{x^{m+1}}{(m+1)!} = e^{-x}\sum_{n=1}^{\infty}n^{p+1}\frac{x^n}{n!}.$$

We now use property (11) of the exponential polynomials in the form

$$\varphi_{p+1}(x) = e^{-x}\sum_{n=0}^{\infty}n^{p+1}\frac{x^n}{n!}$$

and the proof is finished, since the left hand side in the above equation is exactly $\varphi_{p+1}(x)$.

Formula (16) can also be written in the form

(17) $\quad \sum_{k=1}^{n} S(n,k)\frac{x^k}{k} = e^{-x}\sum_{m=0}^{\infty}\left(1^{n-1} + 2^{n-1} + \ldots + m^{n-1}\right)\frac{x^m}{m!}$

where $n \geq 1$.

Using Bernoulli's formula for the sum of powers we can find a different expression for this Stirling transform. Jacob Bernoulli in 1713 found the representation

$$1^p + 2^p + \ldots + (n-1)^p = \frac{1}{p+1}\sum_{j=0}^{p}\binom{p+1}{j}B_j n^{p+1-j}$$

for integers $p \geq 1$. Here $B_j$ are the Bernoulli numbers with $B_1 = \dfrac{-1}{2}$. We shall write this equation in the form

$$(18) \quad 1^p + 2^p + \ldots + (n-1)^p + n^p = n^p + \frac{1}{p+1}\sum_{k=1}^{p+1}\binom{p+1}{k}B_{p+1-k}n^k$$

by adding $n^p$ to both sides and changing the index of summation to $k = p+1-j$. From here

$$\sum_{n=0}^{\infty}(1^p + 2^p + \ldots + n^p)\frac{x^n}{n!} = \sum_{n=0}^{\infty}\frac{n^p x^n}{n!} + \frac{1}{p+1}\sum_{k=1}^{p+1}\binom{p+1}{k}B_{p+1-k}\left\{\sum_{n=0}^{\infty}\frac{n^k x^n}{n!}\right\}$$

$$= e^x\left(\varphi_p(x) + \frac{1}{p+1}\sum_{k=1}^{p+1}\binom{p+1}{k}B_{p+1-k}\varphi_k(x)\right).$$

Setting $p+1 = n$ we have from the above proposition an interesting corollary.

**Corollary 10.** *For every $n \geq 2$*

$$(19) \quad \sum_{k=1}^{n}S(n,k)\frac{x^k}{k} = \varphi_{n-1}(x) + \frac{1}{n}\sum_{k=1}^{n}\binom{n}{k}B_{n-k}\varphi_k(x)$$

in particular, for $x = 1$

$$(20) \quad \sum_{k=1}^{n}S(n,k)\frac{1}{k} = b_{n-1} + \frac{1}{n}\sum_{k=1}^{n}\binom{n}{k}B_{n-k}b_k$$

(see also [1], p. 51, Example 4).

Using the Bernoulli numbers of the second kind $B_n^+$, where $B_n^+ = B_n\,(n \neq 1)$ and $B_1^+ = \dfrac{1}{2}$ we can write (19) in the form

$$(21) \quad \sum_{k=1}^{n}S(n,k)\frac{x^k}{k} = \frac{1}{n}\sum_{k=1}^{n}\binom{n}{k}B_{n-k}^+\varphi_k(x).$$

This shorter formula is convenient for iteration. We divide by $x$ in (21) and integrate to get

(22) $$\sum_{k=1}^{n} S(n,k)\frac{x^k}{k^2} = \frac{1}{n}\sum_{k=1}^{n}\binom{n}{k}B_{n-k}^{+}\int_{0}^{x}\frac{\varphi_k(t)}{t}dt$$

$$= \frac{1}{n}\sum_{k=1}^{n}\binom{n}{k}B_{n-k}^{+}\left\{\sum_{j=1}^{k}S(k,j)\frac{x^j}{j}\right\} = \frac{1}{n}\sum_{k=1}^{n}\binom{n}{k}B_{n-k}^{+}\left\{\frac{1}{k}\sum_{m=1}^{k}\binom{k}{m}B_{k-m}^{+}\varphi_m(x)\right\}.$$

In particular, for $x=1$ we have

(23) $$\sum_{k=1}^{n} S(n,k)\frac{1}{k^2} = \frac{1}{n}\sum_{k=1}^{n}\binom{n}{k}B_{n-k}^{+}\left\{\frac{1}{k}\sum_{m=1}^{k}\binom{k}{m}B_{k-m}^{+}b_m\right\}.$$

**6. Stirling transform of factorials.**

In this section we shall use the geometric polynomials

(24) $$\omega_n(x) = \sum_{k=0}^{n} S(n,k)k!\,x^k$$

( $n = 0,1,...$ ) studied in [5], [7] and [10]. We shall represent the polynomials

(25) $$\sum_{k=1}^{n} S(n,k)(k-1)!\,x^k = \int_{0}^{x}\frac{\omega_n(t)}{t}dt$$

in terms of $\omega_k(x)$. The following proposition is true.

**Proposition 11**. *For every $n > 1$*

(26) $$\sum_{k=1}^{n} S(n,k)(k-1)!\,x^k = (x+1)\omega_{n-1}(x).$$

*When $n=1$ the left hand side equals $x$.*

*Proof.* The short proof is based on the well-known property of Stirling numbers

$$S(n,k) = kS(n-1,k) + S(n-1,k-1).$$

Thus for $n > 1$ we have

$$\sum_{k=1}^{n} S(n,k)(k-1)!x^k = \sum_{k=1}^{n} S(n-1,k)k!x^k + \sum_{k=1}^{n} S(n-1,k-1)(k-1)!x^k$$

$$= \sum_{k=1}^{n-1} S(n-1,k)k!x^k + \sum_{j=0}^{n-1} S(n-1,j)j!x^{j+1} = \sum_{k=1}^{n-1} S(n-1,k)k!x^k + x\sum_{j=1}^{n-1} S(n-1,j)j!x^j$$

$$= (1+x)\omega_{n-1}(x)$$

since $S(n-1,n) = 0$ and $S(n-1,0) = 0$. We also set $k-1 = j$ in the last sum. The case $n = 1$ is obvious.

**Corollary 12**. The *geometric polynomials have the recurrence property* ($n \geq 1$)

(17) $\quad \omega_n(x) = x\omega_{n-1}(x) + (x + x^2)\omega'_{n-1}(x)$.

*Proof.* Follows immediately from the above proposition and equation (14) after differentiation.

This property was first obtained in [10] by a different method.

The numbers

$$\omega_n(1) = \sum_{k=0}^{n} S(n,k)k!$$

are known as the ordered Bell numbers (also Fubini numbers, preferential arrangement numbers). These numbers appear in the next corollary.

**Corollary 13.** *For every* $n = 1, 2, \ldots$ ,

(28) $\quad \displaystyle\sum_{k=1}^{n} S(n,k)(k-1)! = \begin{cases} 2\omega_{n-1}(1) & n > 1 \\ 1 & n = 1 \end{cases}$

(29) $\quad \displaystyle\sum_{k=1}^{n} S(n,k)(k-1)!(-1)^k = \begin{cases} 0 & n > 1 \\ -1 & n = 1 \end{cases}$

For the proof we first set $x = 1$ and then $x = -1$ in (26). Property (29) can be found in Jordan's book [14], on p. 189.

For completeness it is good to mention here that the ordered Bell numbers $\omega_n(1)$ have a series representation. The geometric polynomials have the property

$$\frac{1}{1-x}\omega_n\left(\frac{x}{1-x}\right) = \sum_{k=0}^{\infty} k^n x^k \quad (|x|<1)$$

(see [7] or [10]). For $x = 1/2$ we find

(30) $\quad \omega_n(1) = \sum_{k=0}^{\infty} \frac{k^n}{2^{k+1}}$ .

**Corollary 14.** For every $n \geq 2$

(31) $\quad \sum_{k=2}^{n} S(n,k)(k-2)!(-1)^k = n-1$

For the proof we apply (2) with $\lambda = \mu = 1$ to the function $f(t) = \sum_{n=2}^{\infty}(n-2)!(-1)^n \frac{t^n}{n!}$ .

## 7. Stirling transform of derangement numbers

The derangement numbers are defined by

$$D_n = n!\left(1 - \frac{1}{1!} + \frac{1}{2!} - \frac{1}{3!} + \ldots + \frac{(-1)^n}{n!}\right)$$

$n = 0, 1, \ldots$ . They give the number of permutations of $\{1, 2, \ldots, n\}$ with no fixed points (see Stanley [20], p. 199). This is sequence A000166 in the On-Line Encyclopedia of Integer Sequences. The generating function of the derangement numbers is given by

$$D(x) = \frac{e^{-x}}{1-x} = \sum_{n=0}^{\infty} D_n \frac{x^n}{n!} \quad (|x|<1).$$

We will prove that the (alternating) Stirling transform of these numbers equals the binomial transform of the Bell numbers $b_n$.

**Theorem 15**. *For every $n = 0, 1, \ldots$ we have*

(32) $\quad (-1)^n \sum_{k=0}^{n} S(n,k)(-1)^k D_k = \sum_{k=0}^{n} \binom{n}{k}(-1)^k b_k = 1 + \sum_{j=0}^{n-1}(-1)^{j+1} b_j$ .

*Proof.* We use property (2) with $\lambda = 1, \mu = -1$, replacing $x$ in $D(x)$ by $e^t - 1$

$$D(-(e^t - 1)) = \sum_{n=0}^{\infty} \frac{t^n}{n!} \left\{ \sum_{k=0}^{n} S(n,k)(-1)^k D_k \right\}.$$

At the same time from the definition of $D(x)$

$$D(-(e^t - 1)) = e^{-t} e^{e^t - 1} = e^{-t} \sum_{n=0}^{\infty} b_n \frac{t^n}{n!} = \sum_{n=0}^{\infty} \frac{t^n}{n!} \left\{ \sum_{k=0}^{n} \binom{n}{k} (-1)^{n-k} b_k \right\}$$

$$\sum_{n=0}^{\infty} \frac{t^n}{n!} \left\{ (-1)^n \sum_{k=0}^{n} \binom{n}{k} (-1)^k b_k \right\}$$

using here the generating function of the Bell numbers [6], [7].

$$e^{e^t - 1} = \sum_{n=0}^{\infty} b_n \frac{t^n}{n!}$$

and Euler's series transformation (see [4]) for the last equality. Comparing coefficients we find the first equality of the theorem. The second equality in (32) comes from the following lemma by using the fact that $b_k = \varphi_k(1)$. The lemma presents the alternating binomial transform of the exponential polynomials $\varphi_n(x)$.

**Lemma 16.** *For every $n = 0, 1, ...$ the exponential polynomials $\varphi_k(x)$ satisfy the binomial identity*

(33) $$\sum_{k=0}^{n} \binom{n}{k} (-1)^k \varphi_k(x) = 1 + x \sum_{j=0}^{n-1} (-1)^{j+1} \varphi_j(x).$$

*Proof.* We have the binomial transform

$$\sum_{k=0}^{n} \binom{n}{k} (-1)^k \varphi_{k+1}(x) = (-1)^n x \varphi_n(x)$$

for $n = 0, 1, ...$ . This is entry (11.67) in [4]. Let for $n = 0, 1, ...$

$$A_n(x) = \sum_{k=0}^{n} \binom{n}{k} (-1)^k \varphi_k(x).$$

Then by the properties of the binomial transform (see [4], equation (1.9))

$$A_{n+1}(x) - A_n(x) = \sum_{k=0}^{n}\binom{n}{k}(-1)^{k+1}\varphi_{k+1}(x) = (-1)^{n+1}x\varphi_n(x).$$

From this recurrence where $A_0(x) = 1$ we find

$$A_n(x) - A_0(x) = x\sum_{j=0}^{n-1}(-1)^{j+1}\varphi_j(x), \quad A_n(x) = 1 + x\sum_{j=0}^{n-1}(-1)^{j+1}\varphi_j(x).$$

The lemma is proved.